\documentclass{amsart}
\usepackage{amsmath,amsthm}
\usepackage{amsfonts,amssymb}
\newtheorem{thm}{Theorem}[section]
\newtheorem{cor}[thm]{Corollary}
\newtheorem{lem}[thm]{Lemma}
\newtheorem{prop}[thm]{Proposition}

\newtheorem{defn}[thm]{Definition}

\theoremstyle{remark}

 \def\CD{{\mathcal D}}
 
 \def\CH{{\mathcal H}}

 \def\CO{{\mathcal O}}
 \def\CP{{\mathcal P}}
 \def\CS{{\mathcal S}}

 \def\CC{{\mathbb C}}

 \def\RR{{\mathbb R}}
 \def\ZZ{{\mathbb Z}}

        \def\sign{\operatorname{sign}}
        \def\id{\operatorname{id}}
        \def\rad{\operatorname{rad}}
\newcommand{\wt}{\widetilde}
\newcommand{\wh}{\widehat}

\begin{document}

\title[convolution and maximal function for Dunkl transform]
 {Convolution operator and maximal function for Dunkl transform}

\author{Sundaram Thangavelu and Yuan Xu}
\address{Stat-Math Division\\Indian Statistical Institute\\
8th Mile, Mysore Road\\Bangalore-560 059\\India.}
\email{veluma@math.iisc.ernet.in}
\address{Department of Mathematics\\ University of Oregon\\
    Eugene, Oregon 97403-1222.}\email{yuan@math.uoregon.edu}

\date{\today}
\keywords{Dunkl transforms, reflection invariance, translation operator,
convolution, summability}
\subjclass{42A38, 42B08, 42B15}
\thanks{ST wishes to thank YX for the warm hospitality during his stay
in Eugene. The work of YX was supported in part by the National Science
Foundation under Grant DMS-0201669}

\begin{abstract}
For a family of weight functions, $h_\kappa$, invariant under a finite
reflection group on $\RR^d$, analysis related to the Dunkl transform
is carried out for the weighted $L^p$ spaces. Making use of the
generalized translation operator and the weighted convolution, we study
the summability of the inverse Dunkl transform, including as examples
the Poisson integrals and the Bochner-Riesz means. We also define
a maximal function and use it to prove the almost everywhere convergence.
\end{abstract}

\maketitle

\section{Introduction}
\setcounter{equation}{0}

The classical Fourier transform, initially defined on
$L^1(\RR^d)$, extends to an isometry of $L^2(\RR^d)$ and it
commutes with the rotation group. For a family of weight functions
$h_\kappa$ invariant under a reflection group $G$, there is a
similar isometry of $L^2(\RR^d,h_\kappa^2)$, called Dunkl
transform (\cite{D92}), which enjoys properties similar to those
of the classical Fourier transform. We denote this transform by
$\wh f$ in the following. It is defined by
$$
 \wh f(x) = \int_{\RR^d} E(x,-i y) f(y) h_\kappa^2(y) dy
$$
where the usual character $e^{-i \langle x, y \rangle}$ is replaced
by  $E(x,-iy) = V_\kappa( e^{- i \langle \cdot, y \rangle})(x)$, in
which $V_\kappa$ is a positive linear operator (see the next section).
If the parameter $\kappa =0$ then $h_\kappa(x) \equiv 1$ and
$V_\kappa = \id$, so that $\wh f$ becomes the classical Fourier
transform.

The basic properties of the Dunkl transform have been studied in
\cite{D92,Jeu,Ros,R03} and also in \cite{Ros2,Tri} (see also the
references therein). These studies are mostly for $L^2(\RR^d)$
or for Schwartz class functions.

The purpose of this paper is to develop an $L^p$ theory for
the summability of the inverse Dunkl transform and prove a
maximal inequality that implies almost everywhere convergence.

The classical Fourier transform behaves well with the translation
operator $f \mapsto f(\cdot - y)$, which leaves the Lebesgue
measure on $\RR^d$ invariant. However, the measure
$h_\kappa^2(x)dx$ is no longer invariant under the usual
translation. One ends up with a generalized
translation operator, defined on the Dunkl transform side by
$$
\wh{\tau_y f}(x) = E(y, -i x) \wh f(x), \qquad x \in \RR^d.
$$
An explicit formula for $\tau_y$ is unknown in general. In fact
$\tau_y$ may not even be a positive operator. Consequently even the
boundedness of $\tau_y$ in $L^p(\RR^d;h_\kappa^2)$ becomes a 
challenging problem. At the moment an explicit formula for $\tau_y f$
is known only in two cases: when $f$ is a radial function and when
$G = \ZZ_2^d$. Properties of $\tau_y$ are studied in Section 3.
In particular, the boundedness of the $\tau_y$ for radial functions is 
established. 

For $f, g$ in $L^2(\RR^d; h_\kappa^2)$ their convolution can be defined 
in terms of the translation operator as
$$
 (f *_\kappa g)(x) = \int_{\RR^d} f(y) \tau_x g^\vee(y) h_\kappa^2(y) dy.
$$
Based on a sharp Paley-Wiener theorem we are able to prove that
$ f*_\kappa \phi_\varepsilon $ converges to $f$ in $L^p(\RR^d;h_\kappa^2)$
for certain radial $\phi$ , where $\phi_\varepsilon$ is a
proper dilation of $\phi$. This and other results are given in Section 4.

The convolution $*_\kappa$ can be used to study the summability
of the inverse Dunkl transform. We prove the $L^p$ convergence of
the summability under mild conditions, including as examples Gaussian
means (heat kernel transform), Abel means and the Bochner-Riesz means
for the Dunkl transform in Section 5.

In  Section 6 we define a maximal function and prove that it is
strong type $(p,p)$ for $1 < p \le \infty$ and weak type $(1,1)$.
As usual, the maximal inequality implies almost everywhere
convergence for the summability.

In the case $G=\ZZ_2^d$, the generalized translation operator is 
bounded on $L^p(\RR^d;h_\kappa^2)$. Many of the results proved in the
previous sections hold under conditions that are more relaxed in this
case and the proof is more conventional. This case will be discussed
in Section 7.

The following section is devoted to the preliminaries and background.
The basic properties of the Dunkl transform will also be given.

\section{Preliminaries}
\setcounter{equation}{0}

Let $G$ be a finite reflection group on $\RR^d$ with a fixed positive root
system $R_+$, normalized so that $\langle v, v \rangle =2$ for all $v \in R_+$,
where $\langle x,y \rangle$ denotes the usual Euclidean inner product. For a
nonzero vector $v \in \RR^d$, let $\sigma_v$ denote the reflection with respect
to the hyperplane perpendicular to $v$, $x \sigma_v:= x-2(\langle x,v \rangle
/\|v\|^2) v$, $x \in \RR^d$. Then $G$ is a subgroup of the orthogonal group
generated by the reflections $\{\sigma_v: v \in R_+\}$.

In \cite{D89}, Dunkl defined a family of first order differential-difference
operators, $\CD_i$, that play the role of the usual partial differentiation
for the reflection group structure. Let $\kappa$ be a nonnegative multiplicity
function $v \mapsto \kappa_v$ defined on $R_+$ with the property that
$\kappa_u = \kappa_v$ whenever $\sigma_u$ is conjugate to $\sigma_v$ in $G$;
then $v \mapsto \kappa_v$ is a $G$-invariant function. Dunkl's operators
are defined by
$$
  \CD_i f(x) = \partial_i f(x) + \sum_{v \in R_+} k_v
    \frac{f(x) -  f(x \sigma_v)} {\langle x, v\rangle}
        \langle v,\varepsilon_i\rangle, \qquad 1 \le i \le d,
$$
where $\varepsilon_1, \ldots, \varepsilon_d$ are the standard unit vectors of
$\RR^d$. These operators map $\CP_n^d$ to $\CP_{n-1}^d$, where $\CP_n^d$
is the space of homogeneous polynomials of degree $n$ in $d$ variables. More
importantly, these operators mutually commute; that is, $\CD_i \CD_j =
\CD_j \CD_i$.

Associated with the reflection group and the function $\kappa$ is the
weight function $h_\kappa$ defined by
\begin{equation}\label{eq:2.1}
h_\kappa(x) = \prod_{v \in R_+} |\langle x, v\rangle|^{\kappa_v}, \qquad
   x \in \RR^d.
\end{equation}
This is a positive homogeneous function of degree $\gamma_\kappa:=
\sum_{v \in R_+} \kappa_v$, and it is invariant under the reflection group
$G$. The simplest example is given by the case $G= \ZZ_2^d$ for which
$h_\kappa$ is just the product weight function
\begin{equation*} 
  h_\kappa (x)  = \prod_{i=1}^d |x_i|^{\kappa_i}, \qquad \kappa_i \ge 0.
\end{equation*}
The Dunkl transform is taken with respect to the measure $h_\kappa^2(x) dx$.

There is a linear isomorphism that intertwines the algebra generated by
Dunkl's operators with the algebra of partial differential operators. The
intertwining operator $V_\kappa$ is a linear operator determined uniquely by
$$
V_\kappa \CP_n \subset \CP_n, \quad V_\kappa1=1, \quad
\CD_i V_\kappa = V_\kappa \partial_i, \quad 1\le i\le d.
$$
The explicit formula of $V_\kappa$ is not known in general. For the group
$G = \ZZ_2^d$, it is an integral transform
\begin{equation}\label{eq:2.2}
  V_\kappa f(x) = b_\kappa \int_{[-1,1]^d} f(x_1 t_1,\ldots, x_d t_d)
   \prod_{i=1}^d (1+t_i)(1-t_i^2)^{\kappa_i-1} dt.
\end{equation}
If some $\kappa_i =0$, then the formula holds under the limit relation
$$
 \lim_{\lambda \to 0} b_\lambda \int_{-1}^1 f(t) (1-t)^{\lambda -1} dt
  = [f(1) + f(-1)] /2.
$$
It is known that $V_\kappa$ is a positive operator (\cite{Ros}); that is,
$p \ge 0$ implies $V_\kappa p\ge 0$.

The function $E(x,y): = V_\kappa^{(x)} \left[e^{\langle x,y\rangle}\right]$,
where the superscript means that $V_\kappa$ is applied to the $x$ variable,
plays an important role in the development of the Dunkl transform. Some of
its properties are listed below (\cite{D91}).

\begin{prop} \label{prop:2.1}
For $x, y \in \RR^n$,
\begin{enumerate}
\item $E(x,y) = E(y,x)$;
\medskip

\item $|E(x,y)| \le e^{\|x\| \cdot \|y\|}$, \qquad $x,y \in \CC^n$;
\medskip

\item Let $\nu(z) = z_1^2 +\ldots + z_d^2$, $z_i\in \CC$. For $z,w \in \CC^d$,
\begin{equation*} 
 c_h \int_{\RR^d} E(z,x) E(w,x) h_\kappa^2(x) e^{-\|x\|^2/2} dx =
    e^{(\nu(z) + \nu(w))/2} E(z,w),
\end{equation*}
where $c_h$ is the constant defined by $c_h^{-1} = \int_{\RR^d} h_\kappa^2(x)
 e^{-\|x\|^2/2}dx$.
\end{enumerate}
\end{prop}

In particular, the function
$$
E(x,iy) = V_\kappa^{(x)} \left[ e^{i \langle x,y\rangle} \right],
     \qquad x ,y \in \RR^d,
$$
plays the role of $e^{i \langle x,y\rangle}$ in the ordinary Fourier analysis.
The Dunkl transform is defined in terms of it by
\begin{equation} \label{eq:2.3}
\wh f(y) =  c_h \int_{\RR^d} f(x) E(x,-iy) h_\kappa^2(x)dx.
\end{equation}
If $\kappa = 0$ then $V_\kappa = id$ and the Dunkl transform coincides with
the usual Fourier transform. If $d =1$ and $G=\ZZ_2$, then the Dunkl transform
is related closely to the Hankel transform on the real line. In fact, in
this case,
$$
 E(x,-iy) = \Gamma(\kappa+1/2) (|x y| /2)^{-\kappa+1/2}
   \left[J_{\kappa-1/2}(|xy|)-
       i  \sign (xy) J_{\kappa+1/2}(|xy|)\right],
$$
where $J_\alpha$ denotes the usual Bessel function
\begin{equation} \label{eq:2.4}
 J_\alpha(t) = \left(\frac{t}{2}\right)^\alpha
      \sum_{n=0}^\infty \frac{(-1)^n } {n!\Gamma(n+\alpha+1)}
       \left( \frac{t}{2} \right)^{2n}.
\end{equation}

We list some of the known properties of the Dunkl transform below
(\cite{D92,Jeu}).

\begin{prop} \label{prop:2.2}
\begin{enumerate}
\item For $f\in L^1(\RR^d; h_\kappa^2)$, $\wh f$ is in $C_0(\RR^d)$.

\item When both $ f $ and $\wh f $ are in $L^1(\RR^d; h_\kappa^2) $
we have the inversion formula
$$
  f(x) = \int_{\RR^d} E(ix,y) \wh f(y) h_\kappa^2(y)dy.
$$

\item The Dunkl transform extends to an isometry of 
$L^2(\RR^d; h_\kappa^2)$.

\item For Schwartz class functions
 $f$, $\wh{\CD_j f} (y) = i y_j \wh f (y)$.
\end{enumerate}
\end{prop}

There are two more results that we will need. They require a little more
preparation. First we need the definition of $h$-harmonics.
The $h$-Laplacian is defined by $\Delta_h=\CD_1^2 + \ldots + \CD_d^2$ and
it plays the role similar to that of the ordinary Laplacian. Let $\CP_n^d$
denote the subspace of homogeneous polynomials of degree $n$ in $d$
variables. An $h$-harmonic polynomial $P$ of degree $n$ is a homogeneous
polynomial $P \in \CP_n^d$ such that $\Delta_h P  =0$. Furthermore,
let $\CH_n^d(h_\kappa^2)$ denote the space of $h$-harmonic polynomials
of degree $n$ in $d$ variables and define
$$
 \langle f, g\rangle_\kappa : = a_\kappa \int_{S^{d-1}} f(x) g(x)
    h^2_\kappa(x) d\omega(x),
$$
where $a_\kappa^{-1} = \int_{S^{d-1}} h_\kappa^2(x) d\omega$.
Then $\langle P,Q\rangle_{\kappa} = 0$ for $P \in \CH_n^d(h_\kappa^2)$ and
$Q \in \Pi_{n-1}^d$. The spherical $h$-harmonics are the restriction of
$h$-harmonics to the unit sphere. The standard Hilbert space theory shows
that
$$
L^2(h_\kappa^2) = \sum_{n=0}^\infty\bigoplus \CH_n^d(h_\kappa^2).
$$

Throughout this paper, we fix the value of $\lambda:=\lambda_\kappa$ as
\begin{equation}\label{eq:2.5}
  \lambda := \gamma_\kappa+ \frac{d-2}{2} \qquad \hbox{with} \qquad
   \gamma_\kappa =  \sum_{v\in R_+} \kappa_v.
\end{equation}
Using the spherical-polar coordinates $x = r x'$, where $x' \in S^{d-1}$,
we have
\begin{equation}\label{eq:2.6}
\int_{\RR^d} f(x) h_\kappa^2(x) dx = \int_0^\infty
 \int_{S^{d-1}} f(r x') h_\kappa^2(x') d\omega(x') r^{2\lambda_\kappa+1} dr
\end{equation}
from which it follows that
$$
  c_h^{-1} = \int_{\RR^d} h_\kappa^2(x) e^{-\|x\|^2/2} dx =
    2^{\lambda_\kappa} \Gamma(\lambda_\kappa+1) a_\kappa^{-1}.
$$


The following formula is useful for computing the Dunkl transform of
certain functions (\cite{D92}).

\begin{prop} \label{prop:2.3}
Let $f \in \CH_n^d(h_\kappa^2)$, $y \in \RR^d$ and $\mu > 0$. Then
the function
$$
  g(x) = a_\kappa \int_{S^{d-1}} f(\xi) E(x,-i\mu \xi)
     h_\kappa^2(\xi)d\omega(\xi)
$$
satisfies $\Delta_h g = -\mu^2 g$ and
$$
  g(x) = (-i)^n f\left(\frac{x}{\|x\|}\right)
         \left(\frac{\mu \|x\|}{2}\right)^{-\lambda_\kappa}
           J_{n+\lambda_\kappa} (\mu\|x\|).
$$
\end{prop}

We will also use the Hankel transform $H_\alpha$ defined on the positive
reals $\RR_+$. For $\alpha > -1/2$,
\begin{equation}\label{eq:2.8}
  H_\alpha f(s) := \frac{1}{\Gamma(\alpha+1)}
      \int_0^\infty f(r) \frac{J_\alpha(r s)}{(r s )^\alpha}
         r^{2 \alpha +1} dr.
\end{equation}
The inverse Hankel transform is given by
\begin{equation}\label{eq:2.9}
  f(r) = \frac{1}{\Gamma(\alpha+1)}
  \int_0^\infty H_\alpha f(s) \frac{J_\alpha(r s)}{(r s )^\alpha}
         s^{2 \alpha +1} ds,
\end{equation}
which holds under mild conditions on $f$; for example, it holds if
$f$ is piecewise continuous and of bounded variation in every finite
subinterval of $(0,\infty)$, and $\sqrt{r} f \in L^1(\RR_+)$
(\cite[p. 456]{W}).

\begin{prop} \label{prop:2.4}
If $f(x) = f_0(\|x\|)$, then $\wh f(x) = H_{\lambda_\kappa} f_0(\|x\|).$
\end{prop}

\begin{proof}
This follows immediately from \eqref{eq:2.6} and Proposition \ref{prop:2.3}.
\end{proof}

\section{Generalized translation}
\setcounter{equation}{0}

One of the important tools in the classical Fourier analysis is the convolution
$$
  (f*g)(x) = \int_{\RR^d} f(y)g(x-y) dy,
$$
which depends on the translation $\tau_y: f(x) \mapsto f(x-y)$. There is a
generalized translation for the reflection invariant weight function, which
we study in this section.

\subsection{Basic properties and explicit formulas}

Taking the Fourier transform, we see that the translation
$\tau_y f = f(\cdot -y)$ of $\RR^d$ satisfies
$\wh {\tau_y f}(x) = e^{-i \langle x, y \rangle} \wh f(x)$. Looking at the
Fourier transform side, an analogue of the translation operator for the Dunkl
transform can be defined as follows:

\begin{defn}
Let $y \in \RR^d$ be given. The generalized translation operator
$f \mapsto \tau_y f$ is defined on $ L^2(\RR^d; h_\kappa^2) $ by
the equation
\begin{equation}\label{eq:3.1}
  \wh{\tau_y f}(x) = E(y,-ix) \wh f(x), \qquad x \in \RR^d.
\end{equation}
\end{defn}

Note that the definition makes sense as the Dunkl transform is 
an isometry of $  L^2(\RR^d; h_\kappa^2) $ onto itself and the
function $ E(y,-ix) $ is bounded. When the function $ f $ is in the
Schwartz class the above equation holds pointwise. Otherwise it is to
be interpreted as an equation for $ L^2 $ functions. As an operator on
$ L^2(\RR^d; h_\kappa^2), \tau_y $ is bounded. A priori it is not
at all clear whether the translation operator can be defined for
$ L^p $ functions for $ p $ different from 2. One of the important 
issues is to prove the $ L^p$ boundedness of the translation operator
on the dense subspace of Schwartz class functions. If it can be done then
we can extend the definition to all $ L^p $ functions.

The above definition gives $ \tau_yf $ as an $ L^2 $ function. It is useful
to have a class of functions on which (3.1) holds pointwise. One such class
is given by the subspace
$$ 
A_\kappa(\RR^d) = \{ f \in  L^1(\RR^d; h_\kappa^2): \wh{f} \in
  L^1(\RR^d; h_\kappa^2) \}. 
$$
Note that $A_\kappa(\RR^d)$ is contained in the intersection of 
$L^1(\RR^d;h_\kappa^2)$ and $ L^\infty $ and hence is a subspace of $
L^2(\RR^d; h_\kappa^2)$. For $f\in A_\kappa(\RR^d)$ we have 
\begin{equation}\label{eq:3.2}
\tau_yf(x) = \int_{\RR^d} E(ix,\xi)E(-iy,\xi)\wh{f}(\xi)h_\kappa^2(\xi) d\xi.
\end{equation}

Before stating some properties of the generalized translation operator
let us mention that there is an abstract formula for $\tau_y$ given
in terms of the intertwining operator $V_\kappa$ and its inverse. It takes
the form of (\cite{Tri})
\begin{equation}\label{eq:VV}
  \tau_y f(x) = V_\kappa^{(x)} \otimes V_\kappa^{(y)}
         \left [(V_\kappa^{-1}f)(x-y)\right]
\end{equation}
for $f$ being Schwartz class functions. 
We note that $V_\kappa^{-1}$ satisfies the formula $V_\kappa^{-1} f(x) =
e^{-\langle y, \CD\rangle} f(x) \vert_{y=0}$. The above formula, however,
does not provide much information on $\tau_y f$. The generalized translation
operator has been studied in \cite{Ros,R03,Tri}. In \cite{Tri} the equation
\eqref{eq:VV} is taken as the starting point.

The following proposition collects some of the elementary properties of
this operator which are easy to prove when both $ f $ and $ g $ are from
$ A_\kappa(\RR^d).$

\begin{prop} \label{prop:3.2}
Assume that $ f \in A_\kappa(\RR^d)$ and $ g \in  L^1(\RR^d; h_\kappa^2)$ 
is bounded. Then
\begin{enumerate}
\item
$\displaystyle{\int_{\RR^d} \tau_y f(\xi) g(\xi) h_\kappa^2(\xi) d\xi
     =  \int_{\RR^d} f(\xi) \tau_{-y} g(\xi) h_\kappa^2(\xi) d\xi.}
$
\medskip
\item $\tau_y f(x) = \tau_{-x} f(-y)$.
\end{enumerate} 
\end{prop} 

\begin{proof}
The property (2) follows from the definition since $ E(\lambda x, \xi)
= E(x,\lambda \xi) $ for any $ \lambda \in \CC.$ To prove (1) assume first
that both $ f $ and $ g $ are from $ A_\kappa(\RR^d).$ Then both integrals in
(1) are well defined. From the definition
\begin{align*}
\int_{\RR^d} \tau_y f(\xi) g(\xi) h_\kappa^2(\xi) d\xi
&= \int_{\RR^d} \left( \int_{\RR^d} E(ix,\xi)E(-iy,\xi)\wh{f}(\xi)
h_\kappa^2(\xi) d\xi \right) g(x) h_\kappa^2(x) dx\\ 
&= \int_{\RR^d} \wh{f}(\xi)\wh{g}(-\xi)E(-iy,\xi)h_\kappa^2(\xi) d\xi.
\end{align*} 
We also have 
\begin{align*}
\int_{\RR^d} f(\xi) \tau_{-y} g(\xi) h_\kappa^2(\xi) d\xi
&= \int_{\RR^d} \left( \int_{\RR^d} E(ix,\xi)E(iy,\xi)\wh{g}(\xi)
h_\kappa^2(\xi) d\xi \right) f(x) h_\kappa^2(x) dx\\
&= \int_{\RR^d} \wh{f}(-\xi)\wh{g}(\xi)E(iy,\xi)h_\kappa^2(\xi) d\xi\\
&= \int_{\RR^d} \wh{f}(\xi)\wh{g}(-\xi)E(-iy,\xi)h_\kappa^2(\xi) d\xi.
\end{align*}
This proves (1) when both $ f $ and $ g $ are from $ A_\kappa(\RR^d).$

Suppose now $ f \in A_\kappa(\RR^d) $ but $ g $ is in the intersection of
$ L^1(\RR^d;h_\kappa^2)$ and $ L^\infty$. Note that $ g \in 
L^2(\RR^d;h_\kappa^2) $ and so $ \tau_yg $ is defined as an $ L^2 $ function.
Since $ f $ is in $ L^2(\RR^d;h_\kappa^2) $ and bounded, both integrals are 
finite. The equation 
$$
\int_{\RR^d} f(\xi) \wh{g}(\xi) h_\kappa^2(\xi) d\xi =
\int_{\RR^d} \wh{f}(\xi) g(\xi) h_\kappa^2(\xi) d\xi
$$ 
which is true for Schwartz class functions remains true for $ f, g \in
L^2 (\RR^d;h_\kappa^2) $ as well. Using this we get 
\begin{align*}
\int_{\RR^d} \tau_y f(x) g(x) h_\kappa^2(x) dx
&= \int_{\RR^d} \tau_y f(-x) g(-x) h_\kappa^2(x) dx\\
& = \int_{\RR^d} E(y,-i\xi)\wh{f}(\xi) \wh{g}(-\xi) h_\kappa^2(\xi) d\xi.
\end{align*} 
By the same argument the integral on the right hand side is also given by
the same expression. Hence (1) is proved.
\end{proof}

We need to prove further properties of $ \tau_y.$ In the classical case
the ordinary translation satisfies
$$ \int_{\RR^d} f(x-y) dx = \int_{\RR^d} f(x) dx.$$ Such a property is
true for $\tau_y $ if $ f $ is a Schwartz class function. Indeed
$$ 
\int_{\RR^d} \tau_yf(x)h_\kappa^2(x) dx = \wh{(\tau_yf)}(0)
= \wh{f}(0).
$$ 
Here we have used the fact that $ \tau_y $ takes $ \CS $ into itself. For
$ f \in A_\kappa(\RR^d) $ though $ \tau_yf $ is defined we do not know if
it is integrable. We now address the question whether the above property holds
at least for a subclass of functions.

For this purpose we make use of the following result which gives an
explicit formula for $ \tau_yf $ when $ f $ is radial, see \cite{R03}. We
use the notation $x' = \frac{x}{|x|}$ for non-zero $ x \in \RR^d.$ 

\begin{prop} \label{thm:3.3}
Let $ f \in A_\kappa(\RR^d) $ be radial and let $f(x) = f_0(\|x\|)$. Then
$$
  \tau_y f(x) = V_\kappa \left[ f_0\left(\sqrt{\|x\|^2+\|y\|^2 -
    2 \|x\|\ \|y\| \langle x', \cdot \rangle}\right) \right](y').
$$
\end{prop}

This formula is proved in \cite{R03} for all Schwartz class functions. However
a different proof can be given using expansions in terms of $h-$harmonics. For
that one needs to invert Hankel transforms of $h-$harmonic coefficients of $ f$
of various orders. Once we assume that $ f \in A_\kappa(\RR^d)$ it follows
that all $h-$ harmonic coefficients of $ f $ and their Hankel transforms
are integrable so that inversion is valid. A special case of the above 
theorem is the following formula
\begin{equation} \label{eq:3.3}
  \tau_y q_t(x)  =  e^{-t(\|x\|^2+\|y\|^2)} E(2tx,y)
\end{equation}
where
$$ 
     q_t(x) = (2t)^{-(\gamma+\frac{d}{2})} e^{-t\|x\|^2} 
$$ 
is the so called heat kernel. This formula has already appeared in 
\cite{Ros2}. The other known formula for $\tau_yf$ is the case when 
$G =\ZZ_2^d$. 

\begin{thm} \label{prop:3.4}
Let $ f \in A_\kappa(\RR^d)$ be radial and nonnegative. Then $\tau_yf \geq 0,
\tau_yf \in L^1(\RR^d; h_\kappa^2) $ and 
$$ \int_{\RR^d} \tau_yf(x) h_\kappa^2(x) dx =
 \int_{\RR^d} f(x) h_\kappa^2(x) dx. $$
\end{thm}
\begin{proof} 
As $ f $ is radial, the explicit formula in Proposition \ref{thm:3.3} shows 
that $\tau_yf \geq 0 $ since $ V_\kappa $ is a positive operator. Taking $ g(x)
= e^{-t\|x\|^2} $ and making use of \eqref{eq:3.3} we get
$$ 
\int_{\RR^d} \tau_yf(x) e^{-t\|x\|^2} h_\kappa^2(x) dx 
= \int_{\RR^d} f(x) e^{-t(\|x\|^2+\|y\|^2)}E(\sqrt{2t}x,\sqrt{2t}y)
 h_\kappa^2(x) dx.
$$
As $ |E(x,y)| \leq e^{\|x\|\,\|y\|} $ we can take limit as $ t \rightarrow 0 $
to get 
$$ 
\lim_{t\rightarrow 0} \int_{\RR^d} \tau_yf(x) e^{-t\|x\|^2} h_\kappa^2(x) dx
=\int_{\RR^d} f(x) h_\kappa^2(x) dx.
$$
Since $ \tau_yf \geq 0,$ monotone convergence theorem applied to the integral
on the left completes the proof.
\end{proof}

We would like to relax the condition on $ f $ in the above proposition. In
order to do that we introduce the notion of generalized (Dunkl) convolution.

\begin{defn}
for $ f, g \in L^2(\RR^d; h_\kappa^2) $ we define
$$ 
f*_\kappa g(x) = \int_{\RR^d} f(y)\tau_xg^\vee(y) h_\kappa^2(y) dy$$
where $ g^\vee(y) = g(-y).$
\end{defn}

Note that as $\tau_xg^\vee \in L^2(\RR^d;h_\kappa^2) $ the above convolution 
is well defined. We can also write the definition as
$$ 
f*_\kappa g(x) = \int_{\RR^d} \wh{f}(\xi) \wh{g}(\xi) E(ix,\xi)
 h_\kappa^2(\xi) d\xi.
$$ 
If we assume that $ g $ is also in $ L^1(\RR^d;h_\kappa^2) $ so that
$ \wh{g} $ is bounded, then by Plancherel theorem we obtain
$$
  \|f*_\kappa g\|_{\kappa,2} \leq \|g\|_{\kappa,1} \|f\|_{\kappa,2}.
$$
We are interested in knowing under what conditions on $ g $ the operator
$ f \rightarrow f*_\kappa g $ defined on the Schwartz class can be extended 
to $ L^p (\RR^d;h_\kappa^2)$ as a bounded operator. But now we use the 
$ L^2 $ boundedness of the convolution to prove the following.

\begin{thm} \label{thm:3.4}
Let $ g \in L^1(\RR^d; h_\kappa^2) $ be radial, bounded and nonnegative. Then 
$\tau_yg \geq 0, \tau_yg \in L^1(\RR^d; h_\kappa^2) $ and
$$ 
\int_{\RR^d} \tau_yg(x) h_\kappa^2(x) dx =
 \int_{\RR^d} g(x) h_\kappa^2(x) dx. 
$$
\end{thm}

\begin{proof} 
Let $ q_t $ be the heat kernel defined earlier so that
$ \wh{q_t}(\xi) = e^{-t\|\xi\|^2}.$ By Plancherel theorem
$$
   \|g*_\kappa q_t-g \|_{\kappa,2}^2 = \int_{\RR^d} |\wh{g}(\xi)|^2 
   (1- e^{-t\|\xi\|^2})^2  h_\kappa^2(\xi) d\xi 
$$
which shows that $ g*_\kappa q_t \rightarrow g $ in $ L^2 (\RR^d;h_\kappa^2) $ 
as $ t \rightarrow 0 $.
Since $ \tau_y $ is bounded on $ L^2(\RR^d;h_\kappa^2) $ we have 
$ \tau_y(g*_\kappa q_t) \rightarrow \tau_yg $ in $ L^2 (\RR^d;h_\kappa^2)$ as 
$ t \rightarrow 0 $. By passing to a subsequence if
necessary we can assume that the convergence is also almost everywhere.

Now as $ g $ is radial and nonnegative, the convolution
$$ 
     g*_\kappa q_t(x) = \int_{\RR^d} g(y) \tau_xq_t(y) h_\kappa^2(y) dy 
$$
is also radial and nonnegative. We also note that $ g*_\kappa q_t \in
A_\kappa(\RR^d)$ as $ g $ is both in $ L^1(\RR^d;h_\kappa^2) $ and 
$ L^ 2(\RR^d;h_\kappa^2)$; in fact $g*_\kappa q_t \in L^1(\RR^d,h_\kappa^2)$
as $ q_t \in A_\kappa(\RR^d) $and, by Plancherel theorem and H\"{o}lder's
inequality, $\|\wh{g*_\kappa q_t}\|_{\kappa,1}= \|\wh{g} \cdot 
\wh{q_t}\|_{\kappa,1} \le \|g\|_{\kappa,2} \|q_t\|_{\kappa,2}$. Thus 
by Theorem \ref{thm:3.4} we know that 
$ \tau_y(g*_\kappa q_t)(x) \geq 0$. This gives us 
$$ 
    \lim_{t\rightarrow 0} \tau_y(g*_\kappa q_t)(x) = \tau_yg(x) \geq 0 
$$ 
for almost every $ x $.  Once the nonnegativity of $\tau_yg(x) $ is proved 
it is easy to show that it is integrable. As before
$$ 
    \int_{\RR^d} \tau_yg(x) e^{-t\|x\|^2} h_\kappa^2(x) dx\\
    = \int_{\RR^d} g(x) e^{-t(\|x\|^2+\|y\|^2)}E(\sqrt{2t}x,\sqrt{2t}y)    
     h_\kappa^2(x)dx.
$$
Taking limit as $ t $ goes to $ 0 $ and using monotone convergence theorem
we get 
$$ \int_{\RR^d} \tau_yg(x) h_\kappa^2(x) dx =
 \int_{\RR^d} g(x) h_\kappa^2(x) dx. $$
This completes the proof.
\end{proof} 

There is another way of proving the above result which avoids the intermediate
steps. Assuming that $\int_{\RR^d} g(x) h_\kappa^2(x) dx =1$, our result is
an immediate consequence of Proposition 6.2 in \cite{Tri}.We are thankful to 
the referee for pointing this out. We are now in a position to prove the 
following result. Let $L^p_{\rad}(\RR^d;h_\kappa^2)$ denote the space of 
all radial functions in $L^p(\RR^d;h_\kappa^2)$.

\begin{thm} \label{thm:3.5}
The generalized translation operator $ \tau_y $, initially defined on the
intersection of $ L^1 (\RR^d;h_\kappa^2) $ and $ L^\infty $, can be 
extended to all radial functions in $L^p(\RR^d;h_\kappa^2)$,
$1 \leq p \leq 2 $, and $ \tau_y : L_{\rad}^p(\RR^d, h_\kappa^2) \rightarrow
 L^p(\RR^d; h_\kappa^2) $ is a bounded operator.
\end{thm}

\begin{proof}
For real valued $ f \in L^1(\RR^d;h_\kappa^2)\cap L^\infty $ which  is radial
the inequality $ -|f| \leq f \leq |f| $ together with the nonnegativity
of $ \tau_y $ on radial functions in $ L^1(\RR^d;h_\kappa^2)\cap L^\infty $
shows that $ |\tau_yf(x)| \leq \tau_y|f|(x)$. Hence
$$ 
\int_{\RR^d} |\tau_yf(x)| h_\kappa^2(x) dx \leq 
\int_{\RR^d} |f|(x)h_\kappa^2(x) dx \leq \|f\|_{\kappa,1}.
$$ We also have
$ \|\tau_yf\|_{\kappa,2} \leq \|f\|_{\kappa,2}.$  By interpolation we get,
as $ L^p $ is the interpolation space between $ L^1 $ and $ L^2$, 
$ \|\tau_yf\|_{\kappa,p} \leq \|f\|_{\kappa,p}$ for all $ 1 \leq p \leq 2$
for all $ f \in L_{\rad}^p(\RR^d;h_\kappa^2)$. This proves the theorem.
For the inerpolation theorem used here see \cite{SW}.
\end{proof}

\begin{thm} \label{cor:3.6}
For every $ f \in L^1_{\rad}(\RR^d; h_\kappa^2) $,
$$ 
\int_{\RR^d} \tau_yf(x) h_\kappa^2(x) dx = \int_{\RR^d} f(x) h_\kappa^2(x) dx.
$$
\end{thm}

\begin{proof}
Choose radial functions $ f_n \in A_\kappa(\RR^d) $ so that $ f_n \rightarrow
f $ and $ \tau_yf_n \rightarrow \tau_yf $ in $ L^1(\RR^d;h_\kappa^2).$  Since
$$ \int_{\RR^d} \tau_yf_n(x) g(x)  h_\kappa^2(x) dx = 
 \int_{\RR^d} f_n(x) \tau_{-y}g(x)  h_\kappa^2(x) dx$$ for every $ g \in 
A_\kappa(\RR^d) $ we get, taking limit as $ n $ tends to infinity,
$$
 \int_{\RR^d} \tau_yf(x) g(x)  h_\kappa^2(x) dx =
 \int_{\RR^d} f(x) \tau_{-y}g(x)  h_\kappa^2(x) dx .
$$  
Now take $ g(x) = e^{-t\|x\|^2} $ and take limit as $ t $ goes to $ 0.$. Since
$ \tau_yf \in L^1(\RR^d;h_\kappa^2) $ by dominated convergence theorem we obtain
$$ 
\int_{\RR^d} \tau_yf(x) h_\kappa^2(x) dx =
\int_{\RR^d} f(x) h_\kappa^2(x) dx
$$
for $f\in L^1(\RR^d;h_\kappa^2)$. 
\end{proof}

We remark that it is still an open problem whether $ \tau_y f $ can be defined
for all $ f \in  L^1(\RR^d;h_\kappa^2).$

\subsection{Positivity of $\tau_y$}

As an immediate consequence of the explicit formula for the generalized
translation of radial functions, if $f(x) \in A_\kappa(\RR^d)$ is 
nonnegative, then $\tau_y f(x) \ge 0$ for all $y\in \RR^d$ (\cite{R03}).

One would naturally expect that the generalized translation defines a positive
operator; that is, $\tau_y f (x) \ge 0$ whenever $f(x) \ge 0$. This, however,
turns out not to be the case. For $G=\ZZ_2$, the explicit formula given
in Section 7 shows that $\tau_y$ is not positive in general 
(signed hypergroup, see \cite{R95}).
Below we give an example to show that $\tau_y$ is not positive in a case
where the explicit formula is not available. It depends on a method of
computing generalized translation of simple functions. The explicit
formula (3.2) can be used to define $ \tau_y f $ when $ f $ is a
polynomial.

\begin{lem}
Let $y \in \RR^d$. For $1 \le j \le d$, $\tau_y\{ x_j \} = x_j - y_j$; and
for $1\le j, k \le d$,
$$
 \tau_y \{x_j x_k\} = (x_j -y_j)(x_k - y_k) -
     \kappa \sum_{v \in R_+} \left[V_\kappa (\langle x,y\rangle) -
        V_\kappa (\langle x \sigma_v,y\rangle)\right].
$$
\end{lem}

\begin{proof}
We use \eqref{eq:3.3} and the fact that $\CD_j \tau_y = \tau_y \CD_j$. On
the one hand, since the difference part of $\CD_j$ becomes zero when applied
to radial functions,
$$
\tau_y \CD_j e^{-t\|x\|^2} = -2t \tau_y
   \left(\{\cdot\}_j e^{-t \|\cdot\|^2}\right) (x).
$$
On the other hand it is easy to verify that
$$
 \CD_j \tau_y e^{-t\|x\|^2} =
 \CD_j \left[e^{-t(\|x\|^2 + \|y\|^2)}E(2t x,y) \right]=
     2 t e^{-t(\|x\|^2 + \|y\|^2)}E(2t x,y) (y_j -x_j).
$$
Together, this leads to the equation
\begin{equation} \label{eq:3.4}
 \tau_y \left(x_j e^{-t \|x\|^2}\right)
    = 2 e^{-t(\|x\|^2 + \|y\|^2)}E(2t x,y) (x_j -y_j).
\end{equation}
Taking the limit as $t\to 0$ gives $ \tau_y \{x_j\} = x_j - y_j$.

Next we repeat the above argument, taking \eqref{eq:3.4} as the starting
point. Using the product formula for $\CD_k$ \cite[p. 156]{DX}, a simple
computation gives
\begin{align*}
 \CD_k \tau_y \left(x _j e^{-t \|x\|^2}\right)
 & =  \CD_k \left[(x_j - y_j) e^{-t (\|x\|^2+\|y\|^2)}\right]\\
& =  e^{-t(\|x\|^2+\|y\|^2)}\Big[ -2t (x_j-y_j)(x_k-y_k)E(2tx,y)  \\
& \qquad +  \delta_{k,j} E(2tx,y)
  +  2 \sum_{v\in R_+} \kappa_v \frac{v_k v_j}{\|v\|^2}
    E(2t x \sigma_v,y) \Big].
\end{align*}
On the other hand, computing $\CD_k(x_j e^{-t\|x\|^2})$ leads to
\begin{align*}
 \tau_y \left( \CD_k (x_j e^{-t \|x\|^2})\right)  =
 -2t \tau_y (x_j x_k e^{-t\|x\|^2}) +\tau_y e^{-t\|x\|^2} \Big[ \delta_{k,j}
     + 2 \sum_{v \in R_+} \kappa_v \frac{v_k v_j}{\|v\|^2} \Big].
\end{align*}
Hence, using \eqref{eq:3.3}, the equation $\CD_k \tau_y ( x_j e^{-t \|x\|^2})
=  \tau_y \CD_k ( x_j e^{-t \|x\|^2})$ gives
\begin{align*}
\tau_y \left(x_jx_k e^{-t\|x\|^2} \right) & =
   e^{-t(\|x\|^2+\|y\|^2)}\Big[ (x_j-y_j)(x_k-y_k)E(2tx,y)  \\
& \qquad +  \sum_{v\in R_+} \kappa_v \frac{v_k v_j}{\|v\|^2}
   \frac{ E(2t x,y) - E(2t x \sigma_v,y)}{t}
 \Big].
\end{align*}
Taking the limit as $t \to 0$ gives the formula of $\tau_y \{x_jx_k\}$.
\end{proof}

\begin{prop}
The generalized translation $\tau_y$ is not a positive operator for the
symmetric group $S_d$.
\end{prop}

\begin{proof}
The formula $\tau_y \{x_j x_k\}$ depends on the values of
$V_\kappa x_j$. For symmetric group $S_d$ of $d$ objects, the formula of
$V_\kappa x_j$ is given by (\cite{D98})
$$
  V_\kappa x_j = \frac{1}{d \kappa + 1} (1+x_j + \kappa |x|), \qquad
     |x|  = x_1+\ldots +x_d.
$$
Let $x(j,k)$ denote the transposition of $x_j$ and $x_k$ variables. It follows
that
\begin{align*}
 \tau_y \{ x_j^2 \} & = (x_j - y_j)^2 + \kappa \sum_{k\ne j}
  \left[ V_\kappa (\langle x, y\rangle ) - V_\kappa (\langle x (k,j),
   y\rangle ) \right] \\
 & = (x_j - y_j)^2 + \kappa \sum_{k\ne j}
  \left[ (x_j - x_k) V_\kappa ( y_j - y_k) \right] \\
  & = (x_j - y_j)^2 + \frac{\kappa}{d \kappa+1} \sum_{k\ne j}
  \left[ (x_j - x_k) ( y_j - y_k) \right].
\end{align*}
Choosing $x = (1,0,0,\ldots,0)$ and $y = (0, 2,2, \ldots, 2)$, we see that
$\tau_y (\{\cdot\}_1^2)(x) = - ( (d-2)\kappa+1 )/(d \kappa+1) \le  0$. This
proves the proposition.
\end{proof}

Let us point out that, by \eqref{eq:3.2}, this proposition also shows that
$V_\kappa^{-1}$ is not a positive operator for the symmetric group. In
the case of $\ZZ_2$, an explicit formula of $V_\kappa^{-1}$ is known
(\cite{X02}) which is not positive.

\subsection{Paley-Wiener theorem and the support of $\tau_y$}

In this subsection we prove a sharp Paley-Wiener theorem and study
its consequences. The usual version of Paley-Wiener theorem has been already
proved by de Jeu in his thesis (Leiden,1994). Another type of Paley-Wiener 
theorem has been proved
in \cite{Tri}. Our result is a refined version of the usual Paley-Wiener
which is analogous to an intrinsic version of Paley-Wiener theorem for the
Fourier transform studied by Helgason \cite{Hel}. Recently, a geometric form 
of the Paley-Wiener theorem has been  conjectured and studied in \cite{Jeu2}.

Let us denote by $\{Y_{j,n}: 1 \le j \le \dim \CH_n^d(h_\kappa^2)\}$ an
orthonormal basis of $\CH_n^d(h_\kappa^2)$. First we prove a Paley-Wiener
theorem for the Dunkl transform. 

\begin{thm}
Let $f\in \CS$ and $B$ be a positive number.
Then $f$ is supported in $\{x: \|x\| \le B\}$ if and only if
for every $j$ and $n$, the function
$$
  F_{j,n}(\rho) = \rho^{-n} \int_{S^{d-1}}
    \wh f(\rho x) Y_{j,n} (x) h_\kappa^2(x)d\omega(x)
$$
extends to an entire function of $\rho \in \CC$ satisfying the
estimate
$$
   |F_{j,n} (\rho)| \le c_{j,n} e^{B \|\Im \rho\|}.
$$
\end{thm}

\begin{proof}
By the definition of $\wh f$ and Proposition \ref{prop:2.3},
\begin{align*}
& \int_{S^{d-1}} \wh f(\rho x) Y_{j,n} (x) h_\kappa^2(x)d\omega(x) \\
 & \qquad = c  \int_{\RR^d} \int_{S^{d-1}} E(y, -i \rho x)
    Y_{j,n} (x) h_\kappa^2(x)d\omega(x) f(y) h_\kappa^2(y) dy \\
  & \qquad = c \int_{\RR^d} f(y) Y_{j,n}(y')
    \frac{J_{\lambda_k + n}(\rho \|y\|)}
    {(\rho \|y\|)^{\lambda_k}} h_\kappa^2(y) dy \\
  & \qquad = c \int_0^\infty f_{j,n}(r) \frac{J_{\lambda_k + n}(r \rho)}
    {(r \rho)^{\lambda_k}} r^{2 \lambda_\kappa + n+1} dr,
\end{align*}
where $c$ is a constant and
$$
   f_{j,n}(r) = r^{-n} \int_{S^{d-1}} f(r y') Y_{j,n}(y')
          h_\kappa^2(y')d\omega(y').
$$
Thus, $F_{j,n}$ is the Hankel transform of order $\lambda_\kappa + n$ of
the function $f_{j,n}(r)$. The theorem then follows from the Paley-Wiener
theorem for the Hankel transform (see, for example, \cite{He}).
\end{proof}

\begin{cor}
A function $f \in \CS$ is supported in $\{x:\|x\|\le B\}$
if and only if $\wh f$ extends to an entire function of $\zeta \in \CC^d$
which satisfies
$$
  |\wh f(\zeta)| \le c \,e^{B\|\Im \zeta\|}.
$$
\end{cor}

\begin{proof}
The direct part follows from the fact that $E(x,-i \zeta)$ is entire
and $|E(x,-i \zeta)| \le c\, e^{\|x\| \cdot \|\Im \zeta\|}$. For the
converse we look at
$$
  \int_{S^{d-1}} \wh f(\rho y') Y_{j,n}(y') h_\kappa^2(y') d
\omega(y'), \qquad \rho \in \CC
$$ where $ d\omega $ is the surface measure on $S^{d-1}.$
This is certainly entire and, from the proof of the previous theorem,
has a zero of order $n$ at the origin. Hence,
$$
 \rho^{-n} \int_{S^{d-1}} \wh f(\rho y') Y_{j,n}(y') h_\kappa^2(y') d\omega(y')
$$
is an entire function of exponential type $B$, from which the converse
follows from the theorem.
\end{proof}

\begin{prop} \label{prop:3.10}
Let  $f \in \CS$ be supported in $\{x: \|x\| \le B\}.$ 
Then $\tau_y f$ is supported in $\{x: \|x\| \le B+\|y\| \}$.
\end{prop}

\begin{proof}
Let $g(x) = \tau_y f(x)$. Then $\wh g(\xi) = E(y, -i \xi)\wh f(\xi)$
extends to $\CC^d$ as an entire function of type $B+ \|y\|$.
\end{proof}

This property of $\tau_y$ has appeared in \cite{Tri}. We
note that the explicit formula for  $\tau_y$ shows that the support
set of $\tau_y$ given in Proposition \ref{prop:3.10} is sharp.

An important corollary in this regard is the following result.

\begin{thm} \label{prop:3.11}
If $f \in C_0^\infty(\RR^d)$ is supported in $ \|x\| \leq B $ then
$\|\tau_y f - f \|_{p} \leq c \|y\|(B+\|y\|)^{\frac{N}{p}} $ for 
$1 \le p \le \infty$ where $ N = d+2\gamma_\kappa.$ Consequently,
$ \lim_{y \to 0} \|\tau_y f - f \|_{\kappa,p} = 0.$
\end{thm} 

\begin{proof}
From the definition we have
$$
  \tau_y f (x) - f(x) = \int_{\RR^d} \left(E(y,-i\xi) -1\right)
    E(x,i\xi) \wh f(\xi) h_\kappa^2(\xi) d\xi.
$$
Using mean value theorem and estimates on the derivatives of $ E(x,i\xi)$
we get the estimate
$$
 \|\tau_yf -f\|_\infty \leq c \|y\| \int_{\RR^d} \|\xi\|\, |\wh{f}(\xi)|
 h_\kappa^2(\xi) d\xi.
$$ 
As $ \tau_yf $ is supported in $ \|x\| \leq (B+\|y\|) $ we obtain
$$\|\tau_y f - f \|_{p} \leq c \|y\|(B+\|y\|)^{\frac{N}{p}} $$ which
goes to zero as $ y $ goes to zero.
\end{proof}

\section{The generalized convolution}
\setcounter{equation}{0}

\subsection{Convolution}
Recall that in section 3 we have defined the convolution $ f*_\kappa g $
when $ f, g \in L^2(\RR^d;h_\kappa^2) $ by
$$
  (f *_\kappa g)(x) = \int_{\RR^d} f(y) \tau_xg^\vee(y) h_\kappa^2(y) dy.
$$
This convolution has been considered in \cite{Ros2,Tri}. It satisfies the
following basic properties:
\begin{enumerate}
\item $\wh {f *_\kappa g} = \wh f \cdot \wh g$;

\item  $f *_\kappa g = g *_\kappa f$.
\end{enumerate}
We have also noted that the operator $ f \rightarrow f*_\kappa g$ is bounded on
$ L^2(\RR^d; h_\kappa^2)$ provided $ \wh{g} $ is bounded. We are interested
in knowing under what conditions on $ g $ the operator $ f \rightarrow 
f*_\kappa g $
can be extended to $ L^p $ as a bounded operator. If only the generalized
translation operator can be extended as a
bounded operator on  $L^p(\RR^d;h_\kappa^2)$,
then the convolution will satisfy the usual Young's inequality. At present we
can only say something about convolution with radial functions.

\begin{thm} \label{thm:4.1}
Let $ g $ be a bounded radial function in $ L^1(\RR^d;h_\kappa^2).$ Then
$$ f*_\kappa g(x) = \int_{\RR^d} f(y)\tau_xg^\vee(y) h_\kappa^2(y) dy $$
initially defined on the intersection of $ L^1(\RR^d;h_\kappa^2) $ and 
$ L^2(\RR^d;h_\kappa^2) $ extends to all
$ L^p(\RR^d;h_\kappa^2)$, $1 \leq p \leq \infty $ as a bounded operator.
In particular,
\begin{equation}\label{eq:4.2}
   \|f*_\kappa g\|_{\kappa,p} \leq \|g\|_{\kappa,1} \|f\|_{\kappa,p}.
\end{equation}
\end{thm}

\begin{proof}
For $ g \in L^1(\RR^d;h_\kappa^2) $ which is bounded and radial we have 
$ |\tau_yg| \leq \tau_y|g| $ which shows that
$$ \int_{\RR^d} |\tau_yg(x)|h_\kappa^2(x) dx \leq
\int_{\RR^d} |g(x)|h_\kappa^2(x) dx.$$
Therefore,
$$ \int_{\RR^d} | f*_\kappa g(x)| h_\kappa^2(x) dx \leq \|f\|_{\kappa,1}
 \|g\|_{\kappa,1}.$$
We also have $ \|f*_\kappa g\|_\infty \leq \|f\|_\infty \|g\|_{\kappa,1}.$ 
By interpolation
we obtain $ \|f*_\kappa g\|_{\kappa,p} \leq \|g\|_{\kappa,1} \|f\|_{\kappa,p}$.
\end{proof}

For $\phi \in L^1(\RR^d; h_\kappa^2)$ and $\varepsilon > 0$, we define the
dilation $\phi_\varepsilon$ by
\begin{equation} \label{eq:4.3}
 \phi_\varepsilon(x) = \varepsilon^{- (2\gamma_\kappa+d)}
    \phi(x/\varepsilon).
\end{equation}
A change of variables shows that
$$
  \int_{\RR^d} \phi_\varepsilon (x) h_\kappa^2(x) dx =
      \int_{\RR^d} \phi (x) h_\kappa^2(x) dx, \qquad
      \hbox{for all $\varepsilon > 0$}.
$$

\begin{thm} \label{thm:4.4}
Let $\phi \in L^1(\RR^d; h_\kappa^2)$ be a bounded radial function and 
assume that $c_h \int_{\RR^d}\phi(x)h_\kappa^2(x)dx = 1$. Then for $f \in
L^p(\RR^d;h_\kappa^2)$, $1 \le p < \infty$, and $f \in C_0(\RR^d)$, $p=\infty$,
$$
\lim_{\varepsilon \to 0} \|f *_\kappa \phi_\varepsilon - f \|_{\kappa,p} =0.
$$
\end{thm} 

\begin{proof}
For a given $\eta >0$ we choose $g \in C_0^\infty$ such that
$\|g - f\|_{\kappa,p} < \eta/3$. The triangle inequality and
\eqref{eq:4.2} lead to
$$
  \|f *_\kappa
 \phi_\varepsilon -f \|_{\kappa,p} \le
    \frac{2}{3} \eta +  \|g *_\kappa \phi_\varepsilon -g \|_{\kappa,p}
$$
where we have used $\|g - f\|_{\kappa,p} < \eta/3$.
Since $\phi$ is radial we can choose a radial function $\psi \in C_0^\infty$
such that
$$
\|\phi - \psi\|_{\kappa,1} \le (12 \|g\|_{\kappa,p})^{-1}\eta.
$$
If we let $ a = c_h\int_{\RR^d} \psi(y) h_\kappa^2(y) dy $ then
by the triangle inequality, \eqref{eq:4.2} and \eqref{eq:4.3},
\begin{align*}
 \|g *_\kappa \phi_\varepsilon -g \|_{\kappa,p} & \le
 \|g \|_{\kappa,p} \|\phi - \psi\|_{\kappa,1} +
 \|g *_\kappa \psi_\varepsilon -ag \|_{\kappa,p}+ |a-1| \| g \|_{\kappa,p}\\
 & \le \eta /6 +  \|g *_\kappa \psi_\varepsilon -ag \|_{\kappa,p}
\end{align*}
since $\|g \|_{\kappa,p} \|\phi - \psi\|_{\kappa,1} \le \frac{\eta}{12}$
and
$$ 
|a-1| =\left|c_h\int_{\RR^d} \left(\phi_\varepsilon(x)
  -\psi_\varepsilon(x)\right)
   h_\kappa^2(x) dx \right| \le (12 \|g\|_{\kappa,p})^{-1}\eta.
$$
Thus
$$ 
\|f*_\kappa \phi_\varepsilon - f\|_{\kappa,p} \le \frac{5}{6}\eta + 
\|g *_\kappa \psi_\varepsilon -ag \|_{\kappa,p}.
$$
Hence it  suffices to show that
$\|g *_\kappa \psi_\varepsilon -ag \|_{\kappa,p} \le \eta /6$.

But now $ g \in A_\kappa(\RR^d) $ and so
$$ 
g*_\kappa \phi_\varepsilon(x) =
 \int_{\RR^d} g(y) \tau_x\phi_\varepsilon^\vee(y)
h_\kappa^2(y) dy = \int_{\RR^d}  \tau_{-x}g(y) \phi_\varepsilon(-y)
h_\kappa^2(y) dy .$$ We also know that $ \tau_{-x}g(y)= \tau_{-y}g(x) $
as $ g \in C_0^\infty.$ Therefore,
$$ g*_\kappa \phi_\varepsilon(x) = \int_{\RR^d} \tau_yg(x) \phi_\varepsilon(y) 
h_\kappa^2(y) dy .$$ In view of this
$$ 
 g*_\kappa \psi_\varepsilon(x)-ag(x) =  \int_{\RR^d} 
  \left(\tau_yg(x) -g(x)\right)\psi_\varepsilon(y)h_\kappa^2(y) dy 
$$ 
which gives by Minkowski's integral inequality
$$
\|g*_\kappa \psi_\varepsilon- ag\|_{\kappa,p} \le 
  \int_{\RR^d} \|\tau_yg -g\|_{\kappa,p} 
   \psi_\varepsilon(y)| h_\kappa^2(y) dy. 
$$
If $ g $ is supported in $ \|x\| \le B $ then the estimate in Theorem 
\ref{prop:3.11} gives 
\begin{align*}
\|g*_\kappa \psi_\varepsilon- ag\|_{\kappa,p} & 
   \le c \int_{\RR^d} \|y\| \left(B+\|y\|
  \right)^{\frac{N}{p}}|\psi_\varepsilon(y)| h_\kappa^2(y) dy\\
& \le c \varepsilon \int_{\RR^d} \|y\| 
  \left(B+\|\varepsilon y\|\right)^{\frac{N}{p}} |\psi(y)| h_\kappa^2(y) dy
\end{align*} 
which can be made smaller that $ \frac{\eta}{6} $ by choosing $ \varepsilon$
small. This completes the proof of the theorem.
\end{proof}

The explicit formula in the case of $ G = \ZZ_2^d $ allows us to prove an
analogous result without the assumption that $ \phi $ is radial, see 
Section 7.

\section{Summability of the inverse Dunkl transform}
\setcounter{equation}{0}

Let $\Phi \in L^1(\RR^d;h_\kappa^2)$ be continuous at $ 0$ and assume 
$\Phi(0) =1$. For
$ f \in \CS $ and $\varepsilon > 0$ define
$$
T_\varepsilon f (x) = c_h \int_{\RR^d} \wh f(y) E(ix,y) \Phi(- \varepsilon y)
h_\kappa^2(y)dy.
$$
It is clear that $ T_\varepsilon $ extends to the whole of $ L^2 $ as a
bounded operator which follows from Plancherel theorem.
We study the convergence of $T_\varepsilon f$ as $\varepsilon \to 0$.  Note
that $T_0 f= f$ by the inversion formula for the Dunkl transform. If
$ T_\varepsilon f $ can be extended to all $ f \in L^p(\RR^d;h_\kappa^2)
 $ and if
$T_\varepsilon f \to f$ in $L^p(\RR^d;h_\kappa^2)$,  we say that
the inverse Dunkl transform is $\Phi$-summable.

\begin{prop}
Let $\Phi$ and  $\phi = \wh \Phi$ both belong to $L^1(\RR^d;h_\kappa^2)$.
If $ \Phi $ is radial then
$$
    T_\varepsilon f (x) =  (f *_\kappa \phi_\varepsilon)(x)
$$
for all $ f \in L^2(\RR^d;h_\kappa^2)$ and $\varepsilon > 0$.
\end{prop}

\begin{proof}
Under the hypothesis on $ \Phi $ both $ T_\varepsilon $ and the
operator  taking $ f $ into $  (f *_\kappa \phi_\varepsilon) $ extend
to $ L^2(\RR^d;h_\kappa^2) $ as bounded operators. So it is enough to
verify $ T_\varepsilon f (x) =  (f *_\kappa \phi_\varepsilon)(x)$ for all
$ f $ in the Schwartz class. By the definition of the Dunkl transform,
\begin{align*}
T_\varepsilon f(x) & = c_h \int_{\RR^d} \wh{\tau_{-x} f}(y)
       \Phi(- \varepsilon y) h_\kappa^2(y) dy \\
   & = c_h \int_{\RR^d} \tau_{-x} f (\xi)
        c_h \int_{\RR^d}  \Phi(- \varepsilon y) E(y,-i\xi) h_\kappa^2(y) dy
          h_\kappa^2(\xi) d\xi \\
   & = c_h  \varepsilon^{-(d+2\gamma_\kappa)} \int_{\RR^d} \tau_{-x} f (\xi)
       \wh \Phi(- \varepsilon^{-1}\xi) h_\kappa^2(\xi) d\xi \\
  & = (f *_\kappa \phi_\varepsilon)(x)
\end{align*}
where we have changed variable $\xi \mapsto - \xi$ and used the fact that
$\tau_{-x} f(-\xi) = \tau_\xi f(x)$.
\end{proof}

If the radial function $ \phi $ satisfies the conditions of 
Theorem \ref{thm:4.4}  we obtain the following result.

\begin{thm} \label{thm:5.2}
Let $\Phi(x) \in L^1(\RR^d;h_\kappa^2)$ be radial and assume that 
$\wh \Phi \in L^1(\RR^d; h_\kappa^2)$ is bounded and $\Phi(0)=1$. 
For $f \in L^p(\RR^d;h_\kappa^2)$, $T_\varepsilon f $ 
converges to $f$ in $L^p(\RR^d; h_\kappa^2)$ as $\varepsilon \to 0,$ for
$ 1 \le p <\infty$.
\end{thm}

The following remarks on the above theorem are in order. In general the
convolution $ f*g$ of an $L^p $ function $ f $ with an $ L^1$ function
$ g $ is not defined as the translation operator is not defined for general
$ L^p$ functions even when $ p=1.$ However, when $ g $ satisfies the conditions
of Theorem \ref{thm:4.4} we can define the convolution $ f*g$ by integrating
$ f $ against $ \tau_yg $ which makes sense, see Definition 3.5. It is in 
this sense the above
convolution $ f*\varphi_{\epsilon}$ is to be understood. Then as $ 
f*\varphi_{\varepsilon}$ agrees with $ T_{\varepsilon}f$ on Schwartz functions
and as the convolution operator extends to $L^p $ as a bounded operator our theorem is proved. 
 
We consider several examples. In our first example we take $\Phi$ to
be the Gaussian function, $\Phi(x) =
e^{-\|x\|^2/2}$. By (3) of Proposition \eqref{eq:2.1} with $z = iy$ and
$w =0$, $\wh \Phi(x) = e^{-\|x\|^2/2}$. We choose $\varepsilon =
1/\sqrt{2t}$ and define
$$
 q_t(x) = \Phi_\varepsilon(x) = (2t)^{-(\gamma_k + \frac{d}{2})}
             e^{-\|x\|^2/4t}.
$$
Then $q_t(x)$ satisfies the heat equation for the $h$-Laplacian,
$$
   \Delta_h u(x,t) = \partial_t u(x,t),
$$
where $\Delta_h$ is applied to $x$ variables. For this $\Phi$, our
summability method is just $f *_\kappa q_t$. By \eqref{eq:3.3}, the
generalized translation of $q_t$ is given explicitly by
$$
\tau_y q_t(x) =  (2t)^{-(\gamma_k + \frac{d}{2})} e^{-(\|x\|^2+\|y\|^2)/4t}
  E\left(\frac{x}{\sqrt{2t}}, \frac{y}{\sqrt{2t}}\right)
$$
which is the heat kernel for the solution of the heat equation for
$h$-Laplacian. Then a corollary of Theorem \ref{thm:5.2} gives the
following result in \cite{R00}.

\begin{thm} \label{thm:5.3}
Suppose $f \in L^p(\RR^d;h_\kappa^2)$, $1 \le p < \infty$ or 
$f \in C_0(\RR^d)$, $p=\infty$.
\begin{enumerate}
\item The heat transform
$$
  H_t f(x):= (f*_\kappa q_t)(x) = c_h \int_{\RR^d} f(y) \tau_y q_t(x)
    h_\kappa^2(y)dy, \qquad t > 0,
$$
converges to $f$ in $L^p(\RR^d;h_\kappa^2)$ as $t \to 0$.
\item Define $H_0f(x) = f(x)$. Then the function $H_tf(x)$ solves the
initial value problem
$$
 \Delta_h u(x,t) = \partial_t u(x,t), \qquad u(x,0) = f(x), \qquad
     (x,t) \in \RR^d \times [0,\infty).
$$
\end{enumerate}
\end{thm}

Our second example is the analogue of the Poisson summability, where we
take $\Phi(x) = e^{-\|x\|}$. This case has been studied in \cite{RV}.
In this case, one can compute the Dunkl transform $\wh \Phi$ just as in
the case of the ordinary Fourier transform, namely, using
\begin{equation} \label{exp}
e^{- t} = \frac{1}{\sqrt{\pi}} \int_0^\infty \frac{e^{-u}}{\sqrt{u}}
           e^{- t^2/4u} du,
\end{equation}
and making use of the fact that the transform of Gaussian is itself
(see \cite[p. 6]{SW}). The result is
$$
\wh {e^{-\|x\|}}= c_{d,\kappa}
     \frac{1}{(1+\|x\|^2)^{\gamma_\kappa+\frac{d+1}{2}}}, \qquad
   c_{d,\kappa} = 2^{\gamma_\kappa+\frac{d}{2}}
       \frac{\Gamma(\gamma_\kappa+\frac{d+1}{2})}{\sqrt{\pi}}.
$$
In this case, we define the Poisson kernel as the dilation of $\wh \Phi$,
\begin{equation} \label{poisson1}
 P_\varepsilon(x) : =  c_{d,\kappa}
\frac{\varepsilon}{(\varepsilon^2+\|x\|^2)^{\gamma_\kappa+\frac{d+1}{2}}}.
\end{equation}
Since $\Phi(0) =1$, it is easy to see that $\int P(x,\varepsilon)
h_\kappa^2(x)dx =1$. We have

\begin{thm} \label{thm:5.4}
Suppose $f \in L^p(\RR^d;h_\kappa^2)$, $1 \le p < \infty$, or 
$f \in C_0(\RR^d)$, $p=\infty$. Then the Poisson integral 
$f *_\kappa P_\varepsilon$ converges to $f$ in $L^p(\RR^d;h_\kappa^2)$.
\end{thm}

Again the proof is a corollary of Theorem \ref{thm:5.2}. For $\kappa =0$,
it becomes the Poisson summability for the classical Fourier analysis
on $\RR^d$. We remark that this theorem is already proved in M. Rosler's
habilitation thesis by using a different method. We are thankful to
the referee for pointing this out.

Next we consider the analogue of the Bochner-Riesz means for which
$$
  \Phi(x) = \begin{cases} (1- \|x\|^2)^\delta, & \|x\| \le 1, \\
           0, & \hbox{otherwise} \end{cases}
$$
where $\delta > 0$. As in the case of the ordinary Fourier transform, we take
$\varepsilon = 1/R$ where $R > 0$.  Then the Bochner-Riesz means is defined by
$$
S_R^\delta f (x) = c_h \int_{\|y\|\le R}
  \left(1- \frac{\|y\|^2}{R^2} \right)^\delta \wh f(y) E(ix,y)h_\kappa^2(y)dy.
$$
Recall that we have defined $ \lambda_\kappa = \frac{(d-2)}{2}+\gamma_\kappa$
and $ N = d+2\gamma_\kappa.$

\begin{thm} \label{thm:5.5}
If $f \in L^p(\RR^d; h_\kappa^2)$, $1 \le p < \infty$, or $f \in C_0(\RR^d)$, 
$p=\infty$, and $\delta > (N-1)/2$, then
$$
   \|S_R^\delta f - f\|_{\kappa,p} \to 0, \qquad \hbox{as $R \to \infty$}.
$$
\end{thm}

\begin{proof}
The proof follows as in the case of ordinary Fourier transform
\cite[p. 171]{SW}. From Proposition \ref{prop:2.4} and the
properties of the Bessel function, we have
$$
 \wh \Phi (x) = 2^{\lambda_\kappa} \|x\|^{- \lambda_\kappa - \delta -1}
       J_{\lambda_\kappa + \delta + 1}(\|x\|).
$$
Hence, by $J_\alpha(r) = \CO(r^{-1/2})$, $\wh \Phi \in L^1(\RR^d; h_\kappa^2)$
under the condition $\delta > \lambda_\kappa+1/2 = (N-1)/2$.
\end{proof}

We note that $\lambda_\kappa = (d-2)/2+ \gamma_\kappa$ where $\gamma_\kappa$
is the sum of all (nonnegative) parameters in the weight function. If
all parameters are zero, then $h_\kappa(x) \equiv 1$ and we are back to
the classical Fourier transform, for which the index $(d-1)/2$ is the
critical index for the Bochner-Riesz means. We do not know if the index
$(N-1)/2$ is the critical index for the Bochner-Riesz means of the Dunkl
transforms. 

\section{Maximal function and almost everywhere summability}
\setcounter{equation}{0}

For $ f \in L^2(\RR^d;h_\kappa^2) $ we define the maximal function 
$M_\kappa f$ by
$$
  M_\kappa f(x) =  \sup_{r>0} \frac{1}{ d_\kappa r^{d+2\gamma_\kappa}} 
|f *_\kappa \chi_{B_r}(x)|,
$$
where $ \chi_{B_r}$ is the characteristic function of the ball $ B_r$
of radius $ r $ centered at $ 0 $ and  $d_\kappa = a_\kappa /(d+2\gamma_\kappa)$.
Using \eqref{eq:2.6} we have
$\int_{B_r} h_\kappa^2(y) dy = (a_\kappa /(d+2 \gamma_\kappa))
  r^{d+ 2 \gamma_\kappa}$.
Therefore, we can also write $M_\kappa f(x)$ as
$$
  M_\kappa f(x) = \sup_{r>0} \frac{\left |
      \int_{\RR^d} f(y) \tau_x\chi_{B_r}(y) h_\kappa^2(y) dy \right| }
     {\int_{B_r} h_\kappa^2(y) dy}.
$$
If $ \varphi \in C_0^\infty(\RR^d) $ is a radial function such that
$ \chi_{B_r}(x) \le \varphi(x) $ then from Theorem \ref{thm:3.4} it follows
that $ \tau_y\chi_{B_r}(x) \le \tau_y\varphi(x)$. But $ \tau_y\varphi $ is
bounded; hence  $ \tau_y\chi_{B_r}$ is bounded and compactly supported
so that it belongs to $ L^p(\RR^d;h_\kappa^2)$. This means that the 
maximal function $ M_\kappa f $ is defined for all 
$ f \in L^p(\RR^d;h_\kappa^2)$. We also note that as
$ \tau_y\chi_{B_r} \geq 0 $ we have $ M_\kappa f (x) \leq M_\kappa |f|(x).$ 

\begin{thm} \label{thm:6.1}
The maximal function is bounded on $L^p(\RR^d; h_\kappa^2)$ for $1 < p
\le \infty$; moreover it is of weak type $(1,1)$, that is, for $f \in
L^1(\RR^d;h_\kappa^2)$ and $a >0$,
$$
  \int_{E(a)} h_\kappa^2(x) dx \le \frac{c}{a} \|f\|_{\kappa,1}
$$
where $E(a) = \{ x: M_\kappa f(x) > a\}$ and $c$ is a constant independent
of $a$ and $f$.
\end{thm}

\begin{proof}
Without loss of generality we can assume that $ f \geq 0.$ Let $\sigma =
 d + 2\gamma_\kappa+1$ and define for $j \ge 0$, $B_{r,j} =
\{x: 2^{-j-1} r \le \|x\| \le 2^{-j}r\}$. Then
\begin{align*}
 \chi_{B_{r,j}}(y) & = (2^{-j} r)^\sigma (2^{-j} r)^{-\sigma}
  \chi_{B_{r,j}} (y) \\
 & \le C  (2^{-j} r)^{\sigma-1}
    \frac{2^{-j}r}{((2^{-j}r)^2 + \|y\|^2 )^{\sigma/2}}\chi_{B_{r,j}} (y) \\
 & \le C  (2^{-j} r)^{\sigma-1} P_{2^{-j}r}(y),
\end{align*}
where $P_\varepsilon$ is the Poisson kernel defined in \eqref{poisson1} and
$C $ is a constant independent of $ r $ and $ j.$. Since $\chi_{B_r}$ and 
$P_\varepsilon$ are both bounded, integrable radial functions, it follows 
from Theorem \ref{thm:3.4} that
$$
  \tau_x \chi_{B_{r,j}}(y) \le C (2^{-j} r)^{\sigma-1} \tau_x P_{2^{-j} r}(y).
$$
This shows that for any positive integer $ m $ 
\begin{align*}
\int_{\RR^d} f(y) \sum_{j=0}^m \tau_x \chi_{B_{r,j}}(y)h_\kappa^2(y) dy 
&\leq C \sum_{j=0}^\infty (2^{-j} r)^{\sigma-1} \int_{\RR^d} 
f(y) \tau_x P_{2^{-j} r}(y)h_\kappa^2(y) dy \\
& \leq C\, r^{d+2\gamma_\kappa} \sup_{t>0} f*_\kappa P_t(x).
\end{align*}
As $ \sum_{j=0}^m  \chi_{B_{r,j}}(y) $ converges to $ \chi_{B_r}(y) $ in
$ L^1(\RR^d;h_\kappa^2)$, the boundedness of $ \tau_x $ on 
$ L^1_{\rad} (\RR^d;h_\kappa^2)$ shows that $
\sum_{j=0}^m \tau_x \chi_{B_{r,j}}(y)$ converges to $ \tau_x\chi_{B_r}(y) $ in
$ L^1(\RR^d;h_\kappa^2)$. By passing to a subsequence if necessary we can 
assume that $ 
 \sum_{j=0}^m  \tau_x\chi_{B_{r,j}}(y) $ converges to $ \tau_x\chi_{B_r}(y) $
for almost every $ y.$  Thus all the functions involved are uniformly bounded
by $ \tau_x\chi_{B_r}(y) .$ This shows that $ \sum_{j=0}^m \tau_x 
\chi_{B_{r,j}}(y)$ converges to $ \tau_x\chi_{B_r}(y) $ in 
$ L^{p'}(\RR^d;h_\kappa^2)$ and hence
$$
\lim_{m\to \infty} \int_{\RR^d} f(y)\sum_{j=0}^m 
\tau_x\chi_{B_{r,j}}(y) h_\kappa^2(y) dy = \int_{\RR^d} f(y)\tau_x\chi_{B_r}(y)
 h_\kappa^2(y) dy .
$$
Thus we have proved that 
$$ 
f*_\kappa \chi_{B_r}(x)\leq C r^{d+2\gamma_\kappa} \sup_{t>0} f*_\kappa P_t(x)
$$
which gives the inequality $ M_\kappa f(x) \leq C P^{*}f(x)$, where $P^{*}f(x)=
\sup_{t > 0} f *_\kappa P_t(x)  $ is the maximal function associated to
the Poisson semigroup.

Therefore it is enough to prove the boundedness of $ P^*f.$ Here we follow
a general procedure used in \cite{Stein2}. By looking at the Dunkl transforms
of the Poisson kernel and the heat kernel we infer that
$$ 
f*_\kappa P_t(x) =
\frac{t}{\sqrt{2\pi}}\int_0^\infty (f *_\kappa q_s)(x) e^{-t^2/2s} s^{-3/2} ds,
$$ which implies, as in \cite[p. 49]{Stein2}, that
$$
P^*f(x) \leq C \sup_{t >0} \frac{1}{t} \int_0^t Q_s f(x) ds,
$$
where $ Q_sf(x) = f*_\kappa q_s(x) $ is the heat semigroup. Hence using the Hopf-
Dunford-Schwartz ergodic theorem as in \cite[p. 48]{Stein2},
we get the boundedness of $ P^*f $ on $ L^p(\RR^d;h_\kappa^2)$   for 
$ 1 < p \leq \infty $ and the weak type (1,1).
\end{proof}

The maximal function can be used to study almost everywhere convergence of
$ f*_\kappa \varphi_\epsilon$ as they can be controlled by the Hardy-Littlewood
maximal function $ M_\kappa f$ under some conditions on $\varphi.$
Recall that $ N = d+2\gamma_\kappa$.

\begin{thm} \label{thm:6.2}
Let $ \phi \in A_\kappa(\RR^d) $ be a real valued radial function which
satisfies $ |\phi(x)| \le c (1+\|x\|)^{-N-1}.$ Then 
$$ 
\sup_{\varepsilon >0}|f*_\kappa \phi_\varepsilon(x)| \le c M_\kappa f(x).
$$
Consequently, $ f*_\kappa \phi_\varepsilon (x) \rightarrow f(x) $ for 
almost every $ x $ as $ \varepsilon $ goes to $ 0$ for all $ f $ in 
$ L^p(\RR^d;h_\kappa^2)$, $1 \leq p < \infty$.
\end{thm}

\begin{proof}
We can assume that both $ f $ and $ \phi $ are nonnegative. Writing
$$
\phi_\varepsilon(y) = \sum_{j=-\infty}^\infty \phi_\varepsilon(y) 
\chi_{\varepsilon 2^j \le \|y\| \le \varepsilon 2^{j+1}}(y)
$$
we have
$$
\sum_{j=-m}^m \phi_\varepsilon(y)
\tau_x\chi_{\varepsilon 2^j \le \|y\| \le \varepsilon 2^{j+1}}(y)
\le c \sum_{j=-m}^m (1+\varepsilon 2^j)^{-N-1}
\tau_x\chi_{\varepsilon 2^j \le \|y\| \le \varepsilon 2^{j+1}}(y).
$$ 
This shows that
\begin{align*}
\int_{\RR^d} f(y)  \phi_\varepsilon(y)\sum_{j=-m}^m
\chi_{\varepsilon 2^j \le \|y\| \le \varepsilon 2^{j+1}}(y) h_\kappa^2(y)dy
& \le c \sum_{j=-m}^m (1+\varepsilon 2^j)^{-N-1}(\varepsilon2^j)^N 
M_\kappa f(x)\\
& \leq c M_\kappa f(x).
\end{align*}
Since $ \phi(y) \le c (1+\|y\|)^{-N-1} \le c P_1(y) $ it follows that
$ \tau_x\phi(y) \le c \tau_xP_1(y) $ is bounded.
Arguing as in the previous theorem we can show that the left hand side of the
above inequality converges to $ f*_\kappa \phi_\varepsilon(x).$ Thus we obtain
$$ 
\sup_{\varepsilon>0} |f*_\kappa \phi_\varepsilon(x)| \le c M_\kappa f(x)
$$
from which the proof of almost everywhere convergence follows from the 
standard argument. 
\end{proof}

The above two theorems show that the maximal functions $ M_\kappa f $ and
$ P^*f $ are comparable. As a corollary we obtain almost everywhere
convergence of Bochner-Riesz means.

\begin{cor} \label{cor:6.3}
When $ \delta \geq \frac{N+1}{2} $ the Bochner-Riesz means $ S_R^\delta f(x)
$ converges to $ f(x) $ for almost every $ x $ for all $ f \in 
L^p (\RR^d;h_\kappa^2)$, $1\leq
p <\infty.$
\end{cor}

We expect the corollary to be true for all $ \delta > \frac{(N-1)}{2} $ as in
the case of the Fourier transform. This can be proved if in the above theorem
the hypothesis on $ \phi $ can be relaxed to $ |\phi(x)| \le c (1+\|x\|)^{-N-
\epsilon} $ for some $ \epsilon > 0.$  Since we do not know if $ \tau_y(
(1+\|x\|)^{-N-1}) $ is bounded or not we cannot repeat the proof of the above
theorem.

\section{Product weight function invariant under $\ZZ_2^d$}
\setcounter{equation}{0}

Recall that in the case $G= \ZZ_2^d$, the weight function $h_\kappa$ is 
a product function
\begin{equation} \label{eq:7.1}
  h_\kappa (x)  = \prod_{i=1}^d |x_i|^{\kappa_i}, \qquad \kappa_i \ge 0.
\end{equation}

In this case the explicit formula of the intertwining operator $V_\kappa$ is
known (see \eqref{eq:2.2}) and there is an explicit formula for $\tau_y$. 
The following formula is contained in \cite{R95}, where it is studied under
the context of signed hypergroups. 

\begin{thm} \label{thm:7.1}
For $G = \ZZ_2^d$ and $h_\kappa$ in \eqref{eq:7.1},
$$
 \tau_y f(x) = \tau_{y_1} \cdots \tau_{y_d} f(x),
         \qquad y = (y_1,\ldots,y_d) \in \RR^d,
$$
where for $G= \ZZ_2$ and $h_\kappa(t) = |t|^\kappa$ on $\RR$,
\begin{align} \label{eq:tauZ2}
\tau_s f(t)
=& \frac{1}{2}
      \int_{-1}^1 f\left(\sqrt{t^2+s^2-2 s t u} \right)
  \Big(1+\frac{t-s}{\sqrt{t^2+s^2 - 2 s t u}} \Big) \Phi_\kappa(u)du \\
& + \frac{1}{2} \int_{-1}^1 f\left(-\sqrt{t^2+s^2-2 s t u }\right)
  \Big(1-\frac{t-s}{\sqrt{t^2+s^2 - 2 s t u}} \Big) \Phi_\kappa(u)du,
\notag 
\end{align}
where $\Phi_\kappa(u) = b_\kappa (1+u)(1-u^2)^{\kappa-1}$. Consequently, 
for each $y \in \RR^d$, the generalized translation operator $\tau_y$ for
$\ZZ_2^d$ extends to a  bounded operator on $L^p(\RR^d; h_\kappa^2)$. More 
precisely,
$\|\tau_y f\|_{\kappa,p} \le  3 \|f\|_{\kappa,p}$, $1 \le p \le \infty$.
\end{thm}

Since the generalized translation operator $\tau_y$ extends to a bounded
operator on $L^p(\RR^d;h_\kappa^2)$, many results stated in the previous 
sections can be improved and the proofs can be carried out more conveniently
as in the classical Fourier analysis. In particular, the properties of 
$\tau_y$ given in Proposition \ref{prop:3.2}, Theorem \ref{thm:3.4} and 
Theorem \ref{cor:3.6} all hold under the more relaxed condition of 
$f \in L^1(\RR^d;h_\kappa^2)$. 

The standard proof \cite{Z} can now be used to show that the generalized 
convolution satisfies the following analogous of Young's inequality. 

\begin{prop} \label{prop:7.3}
Let $G = \ZZ_2^d$. Let $p,q,r \ge 1$ and $p^{-1} = q^{-1} +r^{-1}- 1$. 
Assume $f \in L^q(\RR^d, h_\kappa^2)$ and $g \in L^r(\RR^d, h_\kappa^2)$, 
respectively. Then
$$ 
  \|f *_\kappa g \|_{\kappa,p} \le c \|f\|_{\kappa,q} \|g\|_{\kappa,r}
$$
\end{prop}

In the following we give several results that improve the corresponding 
results in the previous sections significantly. We start with an improved
version of Theorem \ref{thm:4.4}. The boundedness of $\tau_y$ allows
us to remove the assumption that $\phi$ is radial.

\begin{thm} \label{thm:7.3}
Let $\phi \in L^1(\RR^d, h_\kappa^2)$ and assume
$\int_{\RR^d}\phi(x)h_\kappa^2(x)dx = 1$. Then for $f \in
L^p(\RR^d;h_\kappa^2)$, $1 \le p < \infty$, or $f \in
C_0(\RR^d)$ if $p = \infty$,
$$
\lim_{\varepsilon \to 0} \|f *_\kappa \phi_\varepsilon - f \|_{\kappa,p} =0,
   \qquad 1 \le p \le \infty.
$$
\end{thm}

\begin{proof}
First we assume that $f \in C_0^\infty(\RR^d)$. By Theorem \ref{prop:3.11},
$\| \tau_y f(x) - f(x) \|_{\kappa,p} \to 0$ as $y \to 0$ for $1 \le p \le
\infty$. In general, for $f \in L^p(\RR^d;h_\kappa^2)$ we write $f =
f_1 + f_2$ where $f_1$ is continuous with compact support and
$\|f_2\|_{\kappa,p} \le \delta$. Then the first term of the inequality
$$
  \|\tau_y f(x) - f(x) \|_{\kappa,p} \le
      \|\tau_y f_1(x) - f_1(x) \|_{\kappa,p} +
       \|\tau_y f_2 (x) - f_2(x) \|_{\kappa,p}
$$
goes to zero as $\varepsilon \mapsto 0$ and the second term is bounded by 
$(1+c)\delta$ as $\|\tau_y f_2\|_{\kappa,p} \le c \|f\|_{\kappa,p}$. 
This proves that $\| \tau_y f(x)-f(x)\|_{\kappa,p} \to 0$ as $y \to 0$. 
We have then
\begin{align*}
& c_h \int_{\RR^d} |f *_\kappa g_\varepsilon(x) - f(x) |^p h_\kappa^2(x) d x \\
& \qquad =
c_h \int_{\RR^d} \left| c_h \int_{\RR^d} (\tau_y f(x)-f(x)) g_\varepsilon(y)
    h_\kappa^2(y) d y  \right|^p h_\kappa^2(x) d x \\
& \qquad \le c_h  \int_{\RR^d} \|\tau_y f- f\|_{\kappa,p}^p
    |g_\varepsilon(x)| h_\kappa^2(x) d x \\
& \qquad = c_h  \int_{\RR^d} \|\tau_{\varepsilon y} f- f\|_{\kappa,p}^p
    |g(x)| h_\kappa^2(x) d x,
\end{align*}
which goes to zero as $\varepsilon \to 0$.
\end{proof}

Our next result is about the boundedness of the spherical means operator. 
As in \cite{MT}, we define the spherical mean operator on
$ A_\kappa(\RR^d) $by
$$
 S_r f(x) : = a_\kappa \int_{S^{d-1}}\tau_{ry} f(x) h_\kappa^2(y) d\omega(y).
$$
The generalized convolution of $ f $ with a radial function can be
expressed in terms of the spherical means $ S_rf.$ In fact,
if $ f \in  A_\kappa(\RR^d) $ and $g(x) = g_0(\|x\|)$
is an integrable radial function then, using the spherical-polar coordinates,
\begin{align*} 
(f *_\kappa g)(x) & = c_h \int_{\RR^d} \tau_y f(x) g(y) h_\kappa^2(y) dy \\
  & = c_h \int_0^\infty r^{2\lambda_\kappa+1} g_0(r)
        \int_{S^{d-1}} \tau_{ry'} f(x)  h_\kappa^2(y') dy' dr \\
  & = \frac{c_h}{a_\kappa} \int_0^\infty
          S_r f(x) g_0(r) r^{2\lambda_\kappa+1} dr.
\end{align*}
We shall make use of this later in this section. Regarding boundedness we have

\begin{thm} \label{thm:4.6}
Let $G = \ZZ_2^d$. For $f \in L^p(\RR^d,h_\kappa^2)$, 
$$
 \|S_r f \|_{\kappa,p} \le c \|f\|_{\kappa,p}, \qquad 1 \le p \le \infty.
$$
Furthermore, $\|S_r f - f\|_{\kappa,p} \mapsto 0$ as $r \to 0+$.
\end{thm}

\begin{proof}
Using H\"older's inequality,
$$
 |S_r f(x)|^p \le a_\kappa \int_{S^{d-1}} |\tau_{ry}f(x)|^p
     h_\kappa^2(y)d\omega(y).
$$
Hence, a simple computation shows that
\begin{align*}
    c_h \int_{\RR^d} |S_r f(x)|^p h_\kappa^2(x) dx
 & \le  c_h \int_{\RR^d} a_\kappa \int_{S^{d-1}} |\tau_{ry}f(x)|^p
     h_\kappa^2(y)d\omega(y) h_\kappa^2(x) dx \\
 & = a_\kappa \int_{S^{d-1}} \|\tau_{ry}f\|_{\kappa,p}^p
     h_\kappa^2(y)d\omega(y) \\
 & \le c \|f\|_{\kappa,p}.
\end{align*}
Furthermore, we have
$$
 \|S_r f-f\|_{\kappa,p}^p \le a_\kappa
  \int_{S^{d-1}} \|\tau_{ry} f - f\|_{\kappa,p}^p h_\kappa^2(y) d\omega(y)
$$
which goes to zero as $ r \to 0$ since $\|\tau_{ry} f - f\|_{\kappa,p} \to 0$.
\end{proof}

We remark that the spherical mean value operator is bounded on $ L^p $ for
any finite reflection group not just for $G = \ZZ_2^d.$ To see this we can
make use of a positive integral representation of the spherical mean operator
proved in \cite{R03}. In fact it easily follows that $ S_r $ is actually
a contraction on $ L^p $ spaces.

The boundedness of $\tau_y f$ in $L^p(\RR^d;h_\kappa^2)$ also allows us
to relax the condition of Theorem \ref{thm:6.2}. 

\begin{thm} \label{thm:7.4}
Set $G = \ZZ_2^d$. Let $\phi(x) = \phi_0(\|x\|)  \in L^1(\RR^d,h_\kappa^2)$
be a radial function. Assume that $\phi_0$ is differentiable,
$\lim_{r \to \infty} \phi_0(r) = 0$ and
$\int_0^\infty r^{2 \lambda_\kappa+2} |\phi_0(r)| dr < \infty$, then
$$
   |(f *_\kappa \phi)(x)| \le c M_\kappa f(x).
$$
In particular, if $\phi \in L^1(\RR^d,h_\kappa^2)$ and 
$c_h\int_{\RR^d} \phi(x) h_\kappa^2(x)dx=1$, then
\begin{enumerate}
\item For $1 \le  p \le \infty$, $f *_\kappa \phi_\varepsilon$
converges to $f$ as $\varepsilon \to 0$ in $L^p(\RR^d;h_\kappa^2)$;
\item For $f \in L^1(\RR^d,h_\kappa^2)$,
$(f *_\kappa \phi_\varepsilon)(x)$  converges to $f(x)$ as
$\varepsilon \to 0$ for almost all $x \in \RR^d$.
\end{enumerate}
\end{thm}

\begin{proof}
By definition of the spherical means $S_t f$, we can also write
$$
  M_\kappa f(x) = \sup_{r > 0}
    \frac{\left| \int_0^r t^{2\lambda_\kappa+1} S_t f(x) dt \right|}
    { \int_0^r t^{2\lambda_\kappa+1} dt }.
$$
Since $|M_\kappa f(x) |  \le c M_\kappa |f|(x)$, we can assume $f(x) \ge 0$.
The assumption on $\phi_0$ shows that
\begin{align*}
 \lim_{r \to \infty} \phi_0(r) \int_0^r S_t f(x) t^{2\lambda_\kappa+1} dt
 & = \lim_{r \to \infty} \phi_0(r) \int_{\RR^d} \tau_yf(x) h_\kappa^2(y) dy\\
 & = \lim_{r \to \infty} \phi_0(r) \int_{\RR^d} f(y) h_\kappa^2(y) dy =0.
\end{align*}
Hence, using the spherical-polar coordinates and integrating by parts, we
get
\begin{align*}
(f*_\kappa \phi)(x)&=\int_0^\infty \phi_0(r) r^{2\lambda_\kappa+1} S_rf(x)dr\\
 & = - \int_0^\infty \int_0^r S_t f(x) t^{2\lambda_\kappa+1} dt \phi'(r) dr,
\end{align*}
which implies that
$$
  |(f*_\kappa \phi)(x)| \le c M_\kappa f(x)
        \int_0^\infty r^{2\lambda_\kappa+2} |\phi_0'(r)| dr,
$$
the boundedness of the last integral proves the maximal inequality.
\end{proof}

As an immediate consequence of the this theorem, the Bochner-Riesz means
converge almost everywhere if $\delta > (N-1)/2$ for $G=\ZZ_2^d$,
which closes the gap left open in Corollary \ref{cor:6.3}. 

We can further enhance Theorem \ref{thm:7.4} by removing the assumption 
that $\phi$ is radial. For this purpose, we make the following simple 
observation about the maximal function: If $f $ is nonnegative then we 
can drop the absolute value sign in the definition of the maximal function, 
even though $\tau_y f$ may not be nonnegative.

\begin{lem} \label{lem:max}
If $f \in L^1(\RR^d, h_\kappa^2)$ is a nonnegative function then
$$
  M_\kappa f(x) = \sup_{r>0} \frac{
      \int_{B_r} \tau_y f(x) h_\kappa^2(y) dy }
     {\int_{B_r} h_\kappa^2(y) dy}.
$$
In particular, if $f$ and $g$ are two nonnegative functions then
$$
M_\kappa f + M_\kappa g =  M_\kappa (f + g).
$$
\end{lem}

\begin{proof}
Since $\tau_y \chi_{B_r}(x)$ is nonnegative, we have
$$
  (f *_\kappa \chi_{B_r}) (x)  =
    \int_{\RR^d} f(y) \tau_y \chi_{B_r}(x) h_\kappa^2(y) dy
$$
is nonnegative if $f$ is nonnegative. Hence we can drop the absolute
value sign in the definition of $M_\kappa f$.
\end{proof}

\begin{thm} \label{thm:almsotZ2}
Set $G = \ZZ_2^d$. Let $\phi \in L^1(\RR^d,h_\kappa^2)$ and let
$\psi(x) = \psi_0(\|x\|)
\in L^1(\RR^d,h_\kappa^2)$ be a nonnegative radial function such that
$|\phi(x)| \le \psi(x)$. Assume that $\psi_0$ is
differentiable, $\lim_{r \to \infty} \psi_0(r) = 0$ and
$\int_0^\infty r^{2 \lambda_\kappa+2} |\psi_0(r)| dr < \infty$.
Then $ \sup_{\epsilon >0}|f*_\kappa \phi_\epsilon(x)|$ is of weak 
type $(1,1)$. In particular, if $\phi \in L^1(\RR^d,h_\kappa^2)$
and $c_h \int_{\RR^d} \phi(x)h_\kappa^2(x)dx=1$, then for 
$f \in L^1(\RR^d,h_\kappa^2)$, $(f *_\kappa \phi_\varepsilon)(x)$
converges to $f(x)$ as $\varepsilon \to 0$ for almost all $x \in \RR^d$.
\end{thm}

\begin{proof}
Since $M_\kappa f(x) \le M_\kappa |f|(x)$, we can assume that $f(x) \ge 0$.
The proof uses the explicit formula for $\tau_y f$. Let us first consider
the case of $d =1$. Since $\psi$ is an even function, $\tau_y \psi$ is given
by the formula 
$$
\tau_y f(x)
= \int_{-1}^1 f\left(\sqrt{x^2+y^2-2 x yt}\right) \Phi_\kappa(t) dt
$$
by \eqref{eq:tauZ2}. Since $(x-y)(1+t) = (x-y t) - (y- x t)$, we have 
$$
  \frac{|x-y|}{\sqrt{x^2+y^2 - 2 x y t}} (1+t) \le 2.
$$
Consequently, by the explicit formula of $\tau_y f$ \eqref{eq:tauZ2}, 
the inequality $|\phi(x)| \le \psi(x)$ implies that ,
$$
 |\tau_y \phi(x)| \le \tau_y \psi (x) + 2 \wt \tau_y \psi(x),
$$
where $\wt \tau_y \psi$ is defined by
$$
 \wt \tau_y \psi(x)=b_\kappa\int_{-1}^1 f\left(\sqrt{x^2+y^2 - 2xy t}\right)
   (1-t^2)^{\kappa -1} dt.
$$
Note that $\wt \tau_y \psi$ differs from $\tau_y \psi$ by a factor of
$(1+t)$ in the weight function. Changing variables $t \mapsto - t$ and
$y \mapsto -y$ in the integrals shows that
\begin{align*}
  \int_{\RR} f(y) \wt \tau_y \psi(x) h_\kappa^2(y) dy
  = \int_{\RR}  F(y) \tau_y \psi(x) h_\kappa^2(y) dy.
\end{align*}
where $F(y) = (f(y)+f(-y)) /2$. Hence, it follows that
\begin{align*}
|(f *_\kappa \phi)(x)|  = \left|
   \int_{\RR} f(y) \tau_y \phi(x) h_\kappa^2(y) dy \right|
 \le  (f*_\kappa \psi)(x) + 2 (F *_\kappa \psi)(x).
\end{align*}
The same consideration can be extended to the case of $\ZZ_2^d$ for $d >1$.
Let $\{e_1, \ldots, e_d\}$ be the standard Euclidean basis. For $\delta_j
= \pm 1$ define $x \delta_j = x -  (1+ \delta_j) x_j e_j$ (that is,
multiplying the $j$-th component of $x$ by $\delta_j$ gives $x \delta_j$).
For $1 \le j \le d $ we define
$$
F_{j_1,\ldots j_k} = 2^{-k}
   \sum_{(\delta_{j_1},\ldots,\delta_{j_k}) \in \ZZ_2^k\}}
       f(x \delta_{j_1} \cdots \delta_{j_k}).
$$
In particular, $F_j(x) = (F(x) + F(x \delta_j))/2$,
$F_{j_1,j_2}(x) = (F(x) + F(x \delta_{j_1}) +F(x \delta_{j_2}) +
F(x \delta_{j_1}\delta_{j_2}) )/4$, and the last sum is over $\ZZ_2^d$,
$F_{1,\ldots,d}(x) = 2^{-d} \sum_{\sigma \in \ZZ_2^d} f(x\sigma)$.
Following the proof in the case of $d=1$ it is not hard to see that
\begin{align*}
|(f *_\kappa \phi)(x)|  \le
   (f*_\kappa \psi)(x) & + 2 \sum_{j=1}^d (F_j *_\kappa \psi)(x) +
     4 \sum_{j_1\ne j_2} (F_{j_1,j_2} *_\kappa \psi)(x)  \\
   &  +  \ldots +  2^d (F_{1,\ldots,d} *_\kappa \psi)(x).
\end{align*}
For $G= \ZZ_2^d$, the explicit formula of $\tau_y$ shows that $M_\kappa f(x)$
is even in each of its variables. Hence, applying the result of the previous
theorem on each of the above terms, we get
\begin{align*}
|(f *_\kappa \phi)(x)|  \le
   M_\kappa f(x) + 2 \sum_{j=1}^d M_\kappa F_j (x) & +
     4 \sum_{j_1\ne j_2} M_\kappa F_{j_1,j_2} (x)\\>
   & + \ldots  + 2^d M_\kappa F_{1,\ldots,d}(x).
\end{align*}

Since all $F_j$ are clearly nonnegative, by Lemma \ref{lem:max}, the last
expression can be written as $M_\kappa H$, where $H$ is the sum of all
functions involved. Consequently, since $\|F_{j_1,\ldots,j_d}\|_{\kappa,1}
\le \|f\|_{\kappa,1}$, it follows that
$$
 \int_{\{x: (f*_\kappa \phi)(x) \ge a\}} h_\kappa^2(y) dy
   \le c \frac{ \|H\|_{\kappa,1}}{a} \le  c_d \frac{ \|f\|_{\kappa,1}}{a}.
$$
Hence, $f *_\kappa \phi$ is of weak type $(1,1)$, from which the almost
everywhere convergence follows as usual.
\end{proof}

We note that we do not know if the inequality
$|(f*_\kappa \phi)(x)| \le c M_\kappa f(x)$ holds in this case, since
we only know $M_\kappa (R(\delta) f)(x) = R(\delta) M_\kappa f(x)
= M_\kappa f(x \delta)$, where $R(\delta) f (x) = f(x \delta)$ for
$\delta \in G$, from which we cannot deduce that
$M_\kappa F_{j_1,\ldots,j_k}(x) \le c M_\kappa f(x)$.

\bigskip\noindent
{\it Acknowledgments.} The authors would like to thank the referee for his/her
careful study of the manuscript and making several useful comments and 
suggestions which have lead to a thorough revision of the paper with more 
rigorous proofs. We also thank him/her for pointing out several references 
which we were not aware of.

\enddocument